
\documentstyle{amsppt} \magnification=\magstep1 \NoBlackBoxes
 \def\vol{\text{Vol
}}   
\def \sgn{\text{sgn}} \topmatter

\title The Fourier transform of order statistics with applications to
Lorentz spaces \endtitle \rightheadtext{The Fourier Transform of Order
Statistics and Lorentz Spaces} \author S. J. Dilworth  and A. L.
Koldobsky \endauthor \address Department of Mathematics, University of
South Carolina, Columbia, SC 29208, U.S.A.  \newline Current address:
Department of Mathematics, Texas A\&M University, College Station,
Texas 77843, U.S.A.  \endaddress \email dilworth\@math.scarolina.edu
\endemail \address Division of Mathematics, Computer Science, and
Statistics, University of Texas at San Antonio, San Antonio, TX 78249,
U.S.A. \endaddress \email koldobsk\@ringer.cs.utsa.edu \endemail

\abstract We present a formula for the Fourier transforms of order
statistics in $\Bbb R^n$ showing that all these Fourier
 transforms are equal up to a constant multiple outside the coordinate
 planes in $\Bbb R^n.$

For $a_1\geq ... \geq a_n\ge0$ and $q>0,$  denote by $\ell_{w,q}^n$ the
$n$-dimensional Lorentz space with the norm $\|(x_1,...,x_n)\| =  (a_1
(x_1^{*})^q +...+ a_n (x_n^{*})^q)^{1/q}$, where
$(x_1^{*},...,x_n^{*})$ is the non-increasing permutation of the
numbers $|x_1|,...,|x_n|.$ We use the above mentioned formula and the
Fourier transform criterion of isometric embeddability of Banach spaces
into $L_q$ \cite{10} to prove that, for $n\geq 3$ and $q\leq 1,$ the
space  $\ell_{w,q}^n$ is isometric to a subspace of $L_q$ if and only
if the numbers $a_1,...,a_n$ form an arithmetic progression. For $q>1,$
all the numbers $a_i$ must be equal so that  $\ell_{w,q}^n = \ell_q^n.$
Consequently, the Lorentz function space $L_{w,q}(0,1)$ is isometric to
a subspace of $L_q$  if and only if  {\it either} $0<q<\infty$ and the
weight $w$ is a constant function  (so that $L_{w,q}= L_q$), {\it or}
$q\le 1$ and $w(t)$ is a decreasing linear function.

Finally, we relate our results to the theory of positive definite
functions.
 \endabstract

 \keywords Fourier transform, isometry, Lorentz
spaces, order statistics \endkeywords

\thanks Both authors were supported in part by the NSF Workshop in
Linear Analysis and Probability held at Texas A\&M University in August
1993 \endthanks \endtopmatter \document \baselineskip=14pt \head 1.
Introduction \endhead

For a vector $x=(x_1,...,x_n)\in \Bbb R^n,$  denote by $x^{*} =
(x_1^{*},...,x_n^{*})$ the non-increasing permutation of the numbers
$|x_1|,...,|x_n|.$ We shall consider the order statistics $x_k^{*}$ as
functions of the variables $x_1,...,x_n.$

For $a_1\geq ... \geq a_n\ge0$ (not all zero) and $q>0,$ the expression
$\|(x_1,...,x_n)\| =  (a_1 (x_1^{*})^q +...+ a_n (x_n^{*})^q)^{1/q}$ is
the norm (or $q$-norm if $q<1$) of an $n$-dimensional weighted Lorentz
space which, as usual, we denote by $\ell_{w,q}^n.$ For an infinite
decreasing sequence $w=\langle a_n \rangle$ of positive weights, for
which $\sum_{n=1}^\infty a_n=\infty$, the Lorentz sequence space
$\ell_{w,q}$ is defined similarly.

Let $I$ denote the interval $(0,1)$ or the interval $(0,\infty)$, and
let $w(t)$ be a positive decreasing function defined on $I$ for which
$\int_0^1 w(t) \, dt=1$ and $\int_0^\infty w(t)\,dt=\infty$ (the latter
condition only for $I=(0,\infty)$).  For $0<q<\infty$, the Lorentz
function space $L_{w,q}(I)$ is the space of equivalence classes of
real-valued measurable functions $f$ on $I$ for which the following
norm (or $q$-norm if $q<1$) $$ \|f\|_{w,q}= \left(\int_I f^*(t)^qw(t)\,
dt\right)^{\frac{1}{q}}$$ is finite, where $f^*$ denotes the
non-increasing rearrangement of $|f|$.  For $0<q\le p <\infty$, the
classical Lorentz spaces $L_{p,q}$ introduced in \cite{15} correspond
to the weights $w(t)=(q/p)t^{\frac{q}{p}-1}$.

Sch\"utt \cite{19} proved  that, if $1\le q<2$, then $L_{w,q}(0,1)$ is
isomorphic to a subspace of $L_q$ if and only if it is a $q$-concave
Banach lattice. For the classical spaces $L_{p,q}$, it follows from
Sch\"utt's result, from Carothers and Dilworth \cite{4}, and from M.
Levy \cite{12,13}, that $L_{p,q}$ is isomorphic to a subspace of $L_q$
if and only if $p=q$ or $0<q\le p <2$. For further isomorphic results
about the subspace structure of the $L_{p,q}$ spaces
 we refer the reader to \cite{3,8,12,13}.

The initial purpose of this work was to check the isometric version of
the above results. We started with the question of whether any
finite-dimensional Lorentz space  $\ell_{w,q}^n$ is isometric to a
subspace of $L_q$? Since we expected a negative answer for $n\geq 3$,
we were going to use the following Fourier transform criterion from
\cite{10}: if $q>0$, where $q$ is not an even integer, and if the
$n$-dimensional quasi-Banach space $E$ is isometric to a subspace of
$L_q$, then the distribution $$\gamma(\xi_1,...,\xi_{n-1}) = {1\over
{(2\pi)^{n-1}c_q}} (\|x\|^q)^{\wedge}(\xi_1,...,\xi_{n-1},1)$$ is a
finite measure on $\Bbb R^{n-1}$ (here $c_q = 2^{q+1} \pi^{1/2} \Gamma
((q+1)/2)/\Gamma(-q/2),$ and the Fourier transform is considered in the
sense of distributions).

Calculating the Fourier transform of the norm of the space
$\ell_{w,q}^n$ we ran into a surprising fact: for every continuous
function $f$ on $\Bbb R$ with a power growth at infinity, all  the
order statistics $\sum_{k=1}^n a_k f(x_ k^{*})$ have equal Fourier
transforms up to a constant multiple outside the coordinate planes in
$\Bbb R^n.$ The same result is true if we consider non-increasing
permutations of $x_1,\dots,x_n$ instead of $|x_1|,\dots,|x_n|$. This
fact has, however, a simple explanation which we present in Section 2.

In Section 3 we apply the formula for the Fourier transform of order
statistics to the Lorentz sequence spaces. We prove that, for $q\leq
1,$ the space  $\ell_{w,q}^n$ is isometric to a subspace of $L_q$ if
and only if the numbers $a_1,\dots,a_n$ form an arithmetic progression.
For $q>1,$ all the numbers $a_i$ must be equal so that  $\ell_{w,q}^n =
\ell_q^n.$

In Section~4 we deduce from the finite-dimensional results that there
are non-trivial isometric embeddings of the function space space
$L_{w,q}(I)$ into  $L_q$  if and only if $I=(0,1)$, $q\le 1$, and
$w(t)$ is a decreasing linear function. As a consequence, we obtain an
interesting family of  rearrangement-invariant renormings of $L_1$
which are at the same time isometric to subspaces of $L_1$.

In Section~5 we show how our results are related to a problem of I. J.
Schoenberg \cite{18} about positive definite functions.

Finally, we wish to mention that the isometries of $L_{p,1}$ into
itself were determined by Carothers and Turett \cite{6}, and that
recently Carothers, Haydon, and Pei-Kee Lin \cite{7} determined the
isometries of $L_{w,q}$ into itself. The methods required to prove
these results are quite different from those used in this paper.

\head 2. The Fourier Transform of Order Statistics \endhead

We start with the following elementary fact.

\proclaim{Lemma 1} For any $a_1, a_2,\dots,a_n\in \Bbb R,$ any function
$f$ on $\Bbb R$, and every $x=(x_1,\dots,x_n)\in \Bbb R^n$, we have
$$\split a_1 f(x_1^{*}) &+\dots+ a_n f(x_n^{*})   \\
 &=\sum_{k=1}^n \Biggl( \sum_{j=1}^k (-1)^{j-1} {{k-1} \choose {j-1}}
 a_{n-k+j} \Biggr) \sum_{i_1<\dots< i_k} f(\max
(|x_{i_1}|,\dots,|x_{i_k}|)). \endsplit \tag{1}
 $$ where the latter sum is taken over all choices of $1\leq
i_1<\dots<i_k\leq n.$ \endproclaim

\demo{Proof} We argue by induction. Assume (without loss of generality)
that $|x_1|\geq |x_k|$ for $k=2,\dots,n$,  and suppose that our
statement is true for the numbers $a_2,\dots,a_n$ and $x_2,\dots,x_n.$
If we add $a_1$ and $x_1$ then the left-hand side of (1) changes by
$a_1 f(|x_1|)$. The additional summands in the right-hand side are as
follows:  $$ a_m \Bigl(\sum_{j=0}^{m-1} (-1)^j {{n-m+j}\choose j}
{{n-1}\choose {n-m+j}}\Bigr) f(|x_1|),  \qquad m=1,\dots,n.$$ For every
$m\geq 2,$ the sum in parentheses is equal to $${{(n-1)\cdot\dots\cdot
m}\over {(n-m)!}} \sum_{j=0}^{m-1} (-1)^j {{m-1}\choose j} =
{{(n-1)\cdot\dots \cdot m}\over {(n-m)!}}((-1)+1)^{m-1} = 0.$$ So the
only non-zero additional summand in the right-hand side is $a_1
f(|x_1|)$, which completes the proof. \qed \enddemo

As usual, we denote by $\Cal S(\Bbb R^n)$ the space of rapidly
decreasing infinitely differentiable functions on $\Bbb R^n,$ and by
$\Cal  S'(\Bbb R^n)$ the space of tempered distributions. For an open
subset $\Omega$ of $\Bbb R^n$, $\Cal D(\Omega)$ denotes the collection
of functions in $\Cal S(\Bbb R^n)$ with compact supports in $\Omega.$
We say that two distributions $f,g\in \Cal  S'(\Bbb R^n)$ are equal on
$\Omega$ if $\langle f, \phi \rangle = \langle g, \phi \rangle$ for
every $\phi\in \Cal D(\Omega)$.

The Fourier transform of any distribution in $\Cal S'(\Bbb R^n)$ of the
form $g(x_{i_1},\dots,x_{i_k})$, where $k<n$, is equal to zero outside
the coordinate planes in $\Bbb R^n.$ Using this fact and Lemma 1 we
immediately get the following result.

\proclaim{Proposition 1} For any continuous function $f$ on $\Bbb R$
with power growth at infinity (i.e., for some $A>0,\rho >0, |f(x)|\leq
A(1+|x|^{\rho})$, for all $x\in \Bbb R$) and any numbers $a_1,
a_2,\dots,a_n\in \Bbb R$ the Fourier transforms of the distributions
$\sum a_k f(x_k^{*})$ and $c f(\max(|x_1|,\dots,|x_n|))$, where $c =
\sum_{k=1}^n (-1)^{k-1} {{n-1} \choose {k-1}} a_{k}$ are equal outside
the coordinate planes in $\Bbb R^n$.

\endproclaim

Lemma 1 also shows that the Fourier transform of order statistics can
easily be calculated if we have a formula for the Fourier transform of
the distributions of the form $f(\max(|x_1|,\dots,|x_n|)).$ Such a
formula was obtained in \cite{11} in connection with some problems
concerning the characterization of measures by potentials. We repeat
this calculation here because the formula is crucial for our further
considerations.

Denote by $G$ the set of vectors $\xi =(\xi_{1},\dots,\xi_{n})\in \Bbb
R^{n}$ such that $\xi_{k}\neq 0$ for $1\leq k\leq n$ and $(\delta ,
\xi) \neq 0$ for every vector $\delta =(\delta_{1},\dots,\delta_{n})$,
with $\delta_{k}=\pm 1$ for $1\leq k\leq n$ (here $(\delta , \xi)$
denotes the usual scalar product in $\Bbb R^n$).

\proclaim{Proposition 2} Let $f$ be an even continuous function on
$\Bbb R^n$ with power growth at infinity and for which the
distribution $u=(f(t)(\sgn (t))^{n-1})^{\wedge}$ is a continuous
function on $\Bbb R\setminus \{0\}.$ Then, for every $\xi\in G,$ we
have $$ \split f(\max(&|x_1|,\dots,|x_n|))^{\wedge}(\xi) \\
&={i^{n-1}\over 2\xi_{1}\cdotp \dots\cdotp\xi_{n}} \sum_{\delta}
\delta_{1}\cdotp \dots\cdotp\delta_{n}
(\delta_{1}\xi_{1}+\dots+\delta_{n}\xi_{n})
u(\delta_{1}\xi_{1}+\dots+\delta_{n}\xi_{n}), \endsplit \tag2$$ where
the sum is taken over all changes of signs.  \endproclaim

\demo{Proof}  Let $\phi \in \Cal S(\Bbb R^{n})$ be a function with a
compact
 support outside the coordinate planes. Then there exists a function
$F\in \Cal S(\Bbb R^{n})$ such that $\partial^{n}F/\partial
x_{1}\dots\partial x_{n}=\widehat {\phi}.$ To see this, note that $\psi
(x)=\phi (x)/(x_{1}\dots x_{n})$ belongs to $\Cal S(\Bbb R^{n})$, and
hence $\widehat {\psi}\in \Cal S(\Bbb R^{n})$. Put $F=i^{n}\widehat
{\psi}$ and now use the connection between the Fourier transform and
differentiation:

$$\partial^{n}F/\partial x_{1}\dots\partial
x_{n}=i^{n}\partial^{n}\widehat {\psi}/\partial x_{1}\dots\partial
x_{n}=(x_{1}\dots x_{n}\psi)^{\wedge}=\widehat {\phi}.$$

Denote by $B_{t}$ the ball $\{x\in \Bbb R^{n}:\|x\|_{\infty}<t\}.$
Recall
 that, for every $\phi \in \Cal S(\Bbb R^{n})$ and for every non-zero
$\xi \in \Bbb R^{n}$, the function $t\rightarrow \widehat{\phi}(t\xi)$
($t\in \Bbb R$)
 is the Fourier transform of the function $y\rightarrow \int_{(\xi ,
x)=y}\phi (x)dx$ ($y\in \Bbb R$): this is the well-known connection
between the Fourier transform and the Radon transform (see e.g.
\cite{9}).  Now we can start the calculation:

$$ \align \langle f(\max(|x_1|,\dots,|x_n|))^{\wedge}, \phi \rangle
 &=\langle f(\max(|x_1|,\dots,|x_n|)), \widehat{\phi} \rangle\\ &=
\int_{R^{n}} f(\max(|x_1|,\dots,|x_n|))\widehat{\phi}(x)\,dx \\
&=\int_{0}^{\infty} f(t)\left(\int_{B_{t}} \widehat{\phi}(x)\,dx
\right)_{t}^{'}\,dt\\ &=\int_{0}^{\infty}
f(t)\left(\int_{-t}^{t}\dots\int_{-t}^{t}
  \widehat{\phi}(x)\,dx\right)_{t}^{'}\,dt \\ &=\int_{0}^{\infty}
f(t)\sum_{\delta} \delta_{1}\cdotp\dots\cdotp\delta_{n}
(F(\delta_{1}t,\dots,\delta_{n}t))_{t}^{'}\,dt, \endalign$$ which
equals $$\int_{0}^{\infty} f(t)\left(\sum_{\delta}
 \delta_{1}\cdotp\dots\cdotp\delta_{n}(\delta_{1}{\partial F\over
 \partial x_{1}}+\dots+\delta_{n}{\partial F\over \partial
x_{n}})(\delta_{1}t,\dots,\delta_{n}t)\right)\,dt. \tag3$$

The function $(\int_{B_{t}} \widehat{\phi}(x)dx)_{t}^{'}$ is even if
$n$ is an odd integer and odd if $n$ is even. Therefore the integral in
(3) is equal to $${1\over 2} \langle f(t)(\sgn (t))^{n-1},
\sum_{\delta}
 \delta_{1}\cdotp\dots\cdotp\delta_{n}\left(\delta_{1}{\partial F\over
\partial x_{1}}+\dots+\delta_{n}{\partial F\over \partial x_{n}}\right)
(\delta_{1}t,\dots,\delta_{n}t)\rangle.\tag {4}$$

Since $\partial F/\partial
 x_{k}=i^{n-1}(x_{k}\phi (x)/(x_{1}\dots x_{n}))^{\wedge}$ for each
 $k$, we can use the property of the Radon transform  which was
mentioned above to rewrite (4) as the following:  $$ {i^{n-1}\over
2}\langle f(t)(\sgn (t))^{n-1},
 \sum_{\delta} \delta_{1}\cdotp\dots\cdotp\delta_{n}
\left(\sum_{j=1}^{n} \int_{(\delta , x)=y} {\delta_{j}x_{j}\phi
(x)\over x_{1}\dots x_{n}} \,dx\right)^{\wedge} \rangle,$$ which equals
$${i^{n-1}\over 2}\langle (f(t)(\sgn (t))^{n-1})^{\wedge}(y),
y\sum_{\delta}
 \delta_{1}\cdotp\dots\cdotp\delta_{n} \int_{(\delta , x)=y} {\phi
 (x)\over x_ {1}\dots x_{n}}\,dx \rangle. $$

The distribution $u=(f(t)(\sgn(t))^{n-1})^{\wedge}$ is a continuous
function on $\Bbb R\setminus \{0\}$, and so, for every $\phi \in \Cal
S(\Bbb R^{n})$ with compact support in $G$, we have $$\langle
f(\|x\|_{\infty})^{\wedge}, \phi \rangle
 ={i^{n-1}\over 2}\int_{\Bbb R} \left(\sum_{\delta} \delta_{1}
\cdotp\dots\cdotp \delta_{n} \int_{(\delta , x)=y} {\phi (x)\over
x_{1}\dots x_{n}}\,dx \right)yu(y)\,dy. \tag {5}$$ The latter integral
converges absolutely because all the functions $$y\rightarrow
\int_{(\delta , x)=y} {\phi (x)\over x_{1}\dots x_{n}}\,dx \qquad(y\in
\Bbb R)$$ belong to $\Cal S(\Bbb R)$ and have compact supports in $\Bbb
R\setminus \{0\}$.  By the Fubini theorem, the integral in the
right-hand side of (5) is equal to $${i^{n-1}\over 2}\int_{\Bbb R^{n}}
\left(\sum_{\delta}
 \delta_{1}\cdotp\dots\cdotp \delta_{n}
{\delta_{1}x_{1}+\dots+\delta_{n}x_{n}\over x_{1}\dots x_{n}}
u(\delta_{1}x_{1}+\dots+\delta_{n}x_{n})\right)\phi(x)\,dx,$$ which
completes the proof. \qed \enddemo

Putting $f(t)=|t|^q$ and $n=3$ in (2) we get the following fact which
was used in \cite{10}  to prove that the space $\ell_{\infty}^3$ is not
isometric to a subspace of $L_q.$

\proclaim{Corollary 1} For every $q>0$ which is not an even integer the
Fourier transform of the function $\max^q (|x_1|,|x_2|,|x_3|)$ is a
continuous homogeneous sign-changing function on the set $G\subset \Bbb
R^3.$  \endproclaim

\demo{Proof} For every  $q>0$ which is not an even integer, we have
$(|t|^q)^{\wedge} (\xi) = c_q |\xi|^{-1-q}$, for all non-zero $\xi \in
\Bbb R$.  Apply the formula (2) with $f(t)=|t|^q$ and $n=3$:  $$ \split
({\text max}^q (|x_1|,|x_2|,|x_3|))^{\wedge} (\xi_1,\xi_2,\xi_3) &=
{-c_q\over {2\xi_1 \xi_2 \xi_3}} \bigl(|\xi_1+\xi_2+\xi_3|^{-q}
\sgn(\xi_1+\xi_2+\xi_3)\\ &\qquad- |\xi_1+\xi_2-\xi_3|^{-q}
 \sgn(\xi_1+\xi_2-\xi_3)\\ &\qquad-|\xi_1-\xi_2+\xi_3|^{-q}
 \sgn(\xi_1-\xi_2+\xi_3)\\ &\qquad+ |\xi_1-\xi_2-\xi_3|^{-q}
\sgn(\xi_1-\xi_2-\xi_3) \bigr)\endsplit$$ for every
$\xi=(\xi_1,\xi_2,\xi_3)\in G.$ Clearly, the Fourier transform has
opposite signs at the points $\xi_1=3$, $\xi_2= \xi_3=1$ and $\xi_1=3$,
$\xi_2=\xi_3=2$.  \qed \enddemo

\head 3. Isometric Embedding of Lorentz  Sequence Spaces into $L_q$
\endhead

We are going to use the Fourier transform criterion mentioned in the
Introduction to characterize the finite- and infinite-dimensional
Lorentz sequence spaces which are isometric to subspaces of $L_q.$  We
formulated this criterion in the
 Introduction as a necessary condition for the existence of an
isometric embedding. As a matter of fact, one can also use this
criterion as a sufficient condition, although there may be some
complications (for details, see \cite{10, Remark 1}).  For Banach
spaces for which
 the distribution $\gamma$ is a finite measure  one obtains the
 following {\it Levy representation} of the norm (named after Paul
Levy):  $$\|x\|^q = \int_{\Bbb R^{n-1}}
|x_1\xi_1+\dots+x_{n-1}\xi_{n-1}+x_n|^q\, d\gamma
(\xi_1,\dots,\xi_{n-1}).$$

It is easy to see \cite{10} that the Levy representation with a finite
measure $\gamma$ implies the existence of an isometric embedding of the
space into $L_q,$ and we shall use this fact in what follows.  In
particular, we shall use the Levy representation for the norm of the
space $\ell_{\infty}^2$  \cite{10}: for each $q<1$ and for all $x,y\in
\Bbb R$, $${\text max}^q(|x|,|y|) = {1\over{2\pi}}\cot(\pi q/2)
\int_\Bbb R |x + y\xi|^q
 {{|\xi-1|^{-q}-|\xi+1|^{-q}}\over\xi}\, d\xi.\tag{6}$$ Note that $
\gamma (\xi)={{|\xi-1|^{-q}-|\xi+1|^{-q}}\over\xi}\, d\xi$ is a finite
measure on $\Bbb R$, which was calculated  (using (2) above) in
\cite{10}
 as the Fourier transform of the function $\max^q(|x|,|y|)$ . If $q>1$
then the function $\gamma$ is still positive, but it is not integrable
around $1$ and $-1.$ Thus, the space $\ell_{\infty}^2$ is isometric to
a subspace of $L_q$ if and only if $q\leq 1.$ For $q=1,$ the Levy
representation is given by the well-known and particularly simple
$$\max(|x|,|y|)=(1/2)(|x+y|+|x-y|).\tag{7}$$

Now we are ready to apply the Fourier transform criterion to Lorentz
spaces.

\proclaim{Theorem 1} (a) Let $0<q\le 1.$ If $n\geq 2$ then the space
$\ell_{w,q}^n$ is isometric to a subspace of $L_q$ if and only if
$$a_1-a_2=a_2-a_3=\dots =a_{n-1}-a_n. \tag{8}$$ (In particular, the
space $\ell_{w,q}^2$ is isometric to a subspace of $L_q$ for every
choice of $a_1, a_2.$) \newline (b) Let $q>1$ and $n\geq 2$. Then the
space $\ell_{w,q}^n$ is isometric to a subspace of $L_q$ if and only if
$a_1=\dots =a_n$ (so that $\ell_{w,q}^n=\ell_q^n).$ \endproclaim

\demo{Proof} Suppose that $n\geq 2$ and  that the space $\ell_{w,q}^n$
is isometric to a subspace of $L_q.$ Let us prove (8) by induction.
The hyperplane in  $\ell_{w,q}^n$ defined by $x_n=0$ is isometric to
the $(n-1)$-dimensional Lorentz space with the weights
$a_1,\dots,a_{n-1}.$ By the induction hypothesis,
$a_1-a_2=a_2-a_3=\dots =a_{n-2}-a_{n-1}=\alpha.$ Denote $a_{n-1}-a_n$
by $\beta.$

For every $k\geq 3$, we have $$\sum_{j=1}^k (-1)^{j-1} {{k-1} \choose
{j-1}} a_{n-k+j} = (-1)^k (\beta - \alpha). \tag{9}$$ To see this, let
$a_{n-k+1}=u$; then $ a_{n-k+j} = u - (j-1)\alpha$ for $j=1,\dots,k-1$,
and $a_n = u - (k-2)\alpha - \beta$. So the sum in (9) is equal to $$u
\sum_{j=1}^k (-1)^{j-1} {{k-1} \choose {j-1}} + (-1)^k \beta +
 \sum_{j=1}^k (-1)^{j} (j-1){{k-1} \choose {j-1}}\alpha - (-1)^k
 \alpha,$$ and it suffices to note that $\sum_{j=1}^k (-1)^{j-1} {{k-1}
\choose {j-1}} = 0$ and $ \sum_{j=1}^k (-1)^{j} (j-1){{k-1} \choose
{j-1}} = 0.$

Apply (1) with $f(t)=|t|^q$ and use (9) to get the following expression
for the norm of the space $\ell_{w,q}^n$:  $$\split
 \|x\|^q = a_1 (x_1^{*})^q +\dots+ a_n (x_n^{*})^q &= a_n (|x_1|^q
 +\dots +|x_n|^q) + \beta \sum_{i<j} {\text max}^q(|x_i|,|x_j|)\\
&\quad+(\beta - \alpha) \sum_{k=3}^n  (-1)^k \sum_{i_1<\dots< i_k}
{\text max}^q (|x_{i_1}|,\dots,|x_{i_k}|). \endsplit \tag{10}$$

Consider the three-dimensional subspace of $\ell_{w,q}^n$ consisting of
the vectors $x\in \Bbb R^n$ for which $x_3=\dots=x_n.$ It follows from
(10) that this subspace is isometric to a space whose norm can be
represented in the form $$\|(x_1,x_2,x_3)\|^q = (\beta - \alpha)
\sum_{k=3}^n (-1)^k {{n-2}\choose {k-2}} {\text max}^q
(|x_1|,|x_2|,|x_3|) + g(x_1,x_2,x_3),$$ where  $g$ is a linear
combination of the functions $|x_i|^q$ and $\max^q(|x_i|,|x_j|)$ and,
therefore, the Fourier transform of $g$ is supported in the coordinate
planes in $\Bbb R^3.$ Note that  $ \sum_{k=3}^n (-1)^k {{n-2}\choose
{k-2}} = -1$, and so $$\|(x_1,x_2,x_3)\|^q = (\alpha - \beta)
 {\text max}^q (|x_1|,|x_2|,|x_3|) + g(x_1,x_2,x_3).$$ This
three-dimensional space will be isometric to a subspace of $L_q.$ By
the Fourier transform criterion discussed above, the distribution $$
\gamma (\xi_1,\xi_2) = (1/((2\pi)^2 c_q))
(\|(x_1,x_2,x_3)\|^q)^{\wedge} (\xi_1, \xi_2, 1)$$ will be a finite
measure in $\Bbb R^2$ (and, in particular, it must be non-negative). On
the other hand, by Corollary 1, the Fourier transform of $ \max^q
(|x_1|,|x_2|,|x_3|)$ is a continuous homogeneous function which changes
its sign on an open set $G\subset \Bbb R^3$, and so $\gamma$ cannot be
non-negative if $\alpha \neq \beta$.  Thus, $\alpha = \beta$, which
concludes the proof of (8).

Conversely, suppose that (8) is satisfied. Then, for $k\geq 3$, all of
the coefficients \linebreak $\sum_{j=1}^k (-1)^{j-1} {{k-1} \choose
{j-1}} a_{n-k+j}$ in (1) are equal to zero , and the formula (1) turns
into $$\align
 \|x\|^q &= a_1 (x_1^{*})^q +\dots+ a_n (x_n^{*})^q \tag{11}\\ &= a_n
(|x_1|^q +\dots+|x_n|^q) + (a_{n-1} - a_n) \sum_{i<j} {\text
max}^q(|x_i|,|x_j|).  \endalign$$

If $q< 1$, then the norm admits the Levy representation with a finite
measure which is the sum of the measures appearing in (6) and of
point-masses $a_n$ at the points $(1,0,\dots,0),\dots, (0,\dots,0,1)$.
So the space $\ell_{w,q}^n$ is isometric to a subspace of $L_q.$ In the
case $q=1$ one easily gets the Levy representation using (7). This
finishes the proof of part (a).

If $q>1$, and if $q$ is not an even integer, then by the remark at the
beginning of Section 3, the Fourier transforms of the functions
$\max^q(|x_i|,|x_j|)$ are not finite measures. So the space
$\ell_{w,q}^n$ can be isometric to a subspace of $L_q$ only if all
these functions are missing in (11), which happens if and only if
$a_{n-1} = a_n.$ In conjunction with (8) this gives  $a_1=\dots =a_n$,
and we have proved part (b). (The case where $q$ is an even integer
will be dealt with below.) \qed \enddemo

The Fourier transform criterion does not work in the case where $q$ is
an even integer. However, statement (b) of Theorem 1 remains true in
this case. The following simple fact shows that, for every $q>1$, the
space $\ell_{w,q}^n$ is {\it smooth} only if all the numbers $a_k$ are
equal. (Thus, for $q>1$, the space $\ell_{w,q}^n$ is isometric to a
subspace of $L_q$ if and only if $\ell_{w,q}^n= \ell_q^n$.)

\proclaim{Lemma 2} For every $q>1$, the space $\ell_{w, q}^n$ is smooth
if and only if $ a_1 = \dots  = a_n.$ \endproclaim

\demo{Proof} Consider the two-dimensional subspace of $\ell_{w, q}^n$
consisting  of the vectors\linebreak $x=(x_1,\dots,x_n)$ for which
$x_1=\dots=x_{n-1}.$ This subspace is isometric to the two-dimensional
Banach space whose norm is given by $$\|(x,y)\|^q =\cases
(a_1+\dots+a_{n-1})|x|^q + a_n|y|^q &\text{for $|x|\geq |y|$}\\
(a_2+\dots+a_n)|x|^q + a_1|y|^q &\text{for $|x|\leq |y|$}.  \endcases$$

Let $a=a_1+\dots+a_{n-1}$, $b=a_n$, $c=a_2+\dots+a_{n}$, and $d=a_1$.
The curves $a|x|^q + b|y|^q = 1$ and $c|x|^q + d|y|^q = 1$ intersect at
the point whose $x$-coordinate equals $(1/(a+b))^{1/q}.$ (Note that
$a+b=c+d$.) If the space is smooth then the derivatives of the
functions $y=({{1-ax^q}\over b})^{1/q}$ and $y=({{1-cx^q}\over
d})^{1/q}$ at the point  $x=(1/(a+b))^{1/q}= (1/(c+d))^{1/q}$ must
agree. This gives $$ \left({{1-ax^q}\over b} \right)^{\frac{1}{q} -
1}\frac{ax^{q-1}}{b}
 = \left({{1-cx^q}\over d} \right)^{\frac{1}{q} - 1}
\frac{cx^{q-1}}{d},$$ and, therefore, $a/b=c/d$.  This means that
$a_1=a_n$, and we are done. \qed \enddemo

Because a decreasing arithmetic progression of positive numbers is
either constant or finite, we immediately deduce from Theorem~1 the
following corollary.  \proclaim{Corollary 2} Let $0<q<\infty$.  The
Lorentz sequence space $\ell_{w,q}$ is isometric to a subspace of $L_q$
if and only if $w$ is a constant sequence. \endproclaim \subsubhead{A
geometrical argument} \endsubsubhead For $q=1,$ we should like to
sketch a simple geometrical proof of Theorem 1.  First we prove by
induction that $\ell_{w,1}^n$ is isometric to a subspace of $L_1$ only
if the weights $a_k$ satisfy (8).  Let $E$ be a Lorentz space of
dimension $n+1$ with weights $a_1,\dots,a_{n+1}$. If $E$ is isometric
to a subspace of $L_1$ then (8) holds by hypothesis, and so it suffices
to prove that $a_{n-1}-a_n=a_n-a_{n+1}$ to complete the induction. In
particular, we may assume that $a_{n-1}$, $a_n$ and $a_{n+1}$ are not
all equal.  Let $B^*$ denote the unit ball of $E^*$.  Since the unit
ball of $E$ is a polytope, and since $E$ is isometric to a subspace of
$L_1$, it follows (see e.g. \cite{1}) that $B^*$ is a {\it zonotope}
(that is, a Minkowski sum of line segments). By \cite{1, Theorem 3.3},
all of the two-dimensional faces of $B^*$ are centrally symmetric. It
is easily seen that
 the extreme points of $B^*$ are all the sign-changed permutations of
 the vector $ \bold a=(a_1,a_2,\dots,a_{n+1})$.  In particular, one of
the two-dimensional faces of $B^*$ has as its vertices all of the
vectors obtained by permuting the last three coordinates of $\bold a$.
If $a_{n-1}=a_n$ or if $a_n=a_{n+1}$, then this face is triangular,
which contradicts the central symmetry requirement.  If
$a_{n-1}>a_n>a_{n+1}$, then the face is {\it hexagonal}, and the
symmetry condition forces $a_{n-1}-a_n=a_n-a_{n+1}$ as required.  Thus
(8) is a necessary condition.  To show that (8) is also sufficient, one
can check that if (8) is satisfied then $B^*$ has four kinds of
two-dimensional faces: two classes of quadrilateral faces, one class of
octagonal faces, and one class of hexagonal faces like the face
described above. The first three kinds of faces are automatically
centrally symmetric without any condition on the weights, while (8)
guarantees that the hexagonal faces are also centrally symmetric.  So,
by \cite{1, Theorem 3.3} once again, if (8) is satisfied, then $B^*$ is
a zonotope, and hence $E$ is isometric to a subspace of $L_1$. \qed
\head 4. Isometric Embedding of Lorentz Function Spaces into $L_q$
\endhead For $0\le \alpha \le 2$, let
$w_\alpha(t)=1+\frac{\alpha}{2}-\alpha t$.  Observe that the
$w_\alpha$'s are precisely the decreasing linear weights on $[0,1]$
which satisfy $\int_0^1 w(t) \, dt=1$.  \proclaim{Theorem 2}(a) Let
$0<q\le 1$. Then $L_{w,q}(0,1)$ is isometric to a subspace of $L_q$ if
and only if $w=w_\alpha$ for some $\alpha \in [0,2]$.  \newline (b) Let
$q>1$. Then $L_{w,q}(0,1)$ is isometric to a subspace of $L_q$ if and
only if $w=w_0$ (so that $L_{w,q}=L_q$).  \endproclaim \demo{Proof} (a)
First we prove necessity. For $n\ge 1$ and for $1\le k\le n$, let
$e_{n,k}=\chi_{[\frac{k-1}{n},\frac{k}{n})}$, and let $X_n$ be the
linear span of the $e_{n,k}$'s in $L_{w,q}$. Clearly, $X_n$ is
isometric to a Lorentz space with weights $a_{n,k}=\int_{(k-1)/n}^{k/n}
w(t)\, dt$.  By Theorem 1, the sequence $\langle
a_{n,k}\rangle_{k=1}^n$ forms a decreasing arithmetic progression for
each $n$.  This clearly forces $w$ to be a decreasing linear function,
so that $w=w_\alpha$ for some $\alpha \in [0,2]$. For sufficiency,
observe that if $w=w_\alpha$, then by Theorem 1 each $X_n$ is isometric
to a subspace of $L_q$. Since $\cup_{n\ge 1}X_n$ is dense in $L_{w,q}$,
it follows from the fact that the $L_q$ spaces are stable under the
operation of taking ultrapowers that $L_{w,q}$ is isometric to a
subspace of $L_q$ (e.g.  \cite{14, pp. 121-122}). \newline (b) In this
case, by Theorem 1, $\langle a_{n,k} \rangle_{k=1}^n$ is a constant
sequence for each $n$, which forces $w$ to be constant.  Sufficiency is
obvious in this case. \qed \enddemo For the function spaces
$L_{w,q}(0,\infty)$ the situation is different.  Arguing as above and
as in Corollary~2 one immediately obtains the following corollary.
\proclaim{Corollary 3} Let $0<q<\infty$. Then $L_{w,q}(0,\infty)$ is
isometric to a subspace of $L_q$ if and only if $w(t)$ is a constant
function.  \endproclaim \remark{Remarks}1. Clearly, the
$L_{w_\alpha,1}$ spaces are all
 isomorphic to $L_1$ (in fact, the $L_{w_\alpha,1}$ and $L_1$ norms are
equivalent).  On the other hand, the proof of a result of Carothers,
Dilworth and Trautman
 \cite{5, Theorem~2.3} shows that there is no isometry from
 $L_{w_\alpha,1}$ onto $L_{w_\beta,1}$ if $0\le \alpha<\beta \le 2$.
 \newline 2. For $\alpha>0$, there is no isometry from $L_1$ {\it into}
$L_{w_\alpha,1}$. In fact, by \cite{5, Lemma~2.1} there is not even an
isometry from $\ell_1^2$ into $L_{w_\alpha,1}$. \newline 3. Observe
that $L_{w_2,1}$ has a certain universality property among the
$L_{w_\alpha,1}$ spaces.  To be precise, if $0<\alpha\le 2$, then
$L_{w_\alpha,1}$ is isometric by a dilation to the one-complemented
sublattice of $L_{w_2,1}$ which consists of all functions that are
supported in the interval $[0,\frac{2\alpha}{\alpha+2}]$.  Moreover,
$L_{w_2,1}$ contains a {\it continuum} of isometrically distinct
one-complemented (by conditional expectation) two-dimensional
subspaces: namely, the Lorentz spaces corresponding to $a_1=1$ and $a_2
\in [1/3,1)$ (whose unit balls are {\it octagonal}).  This is in
contrast with $L_1$, which has, for each positive integer $n$, a unique
one-complemented  $n$-dimensional subspace up to isometry, namely
$\ell_1^n$.  \newline 4. Let us recall some of what is known about the
geometry of the $L_{w,1}$ spaces. For every strictly decreasing $w$,
the extreme points of the unit ball of $L_{w,1}$ are all the functions
of the form $\varepsilon \chi_A /(\int_0^{|A|}w(t)\, dt)$, where $A$ is
any subset of $[0,1]$ with positive Lebesgue measure $|A|$, and
$\varepsilon$ is any $\pm1$-valued measurable function (\cite{5, Prop.
2.2}). Moreover, while $L_{w,1}$ is not strictly convex, it is
nevertheless true that every element on the unit sphere of $L_{w,1}$ is
the barycenter of a {\it unique} Borel probability measure supported on
the extreme points of the ball (\cite{5, Theorem 3.5}).  Finally,
Sedaev \cite{20} proved that if $w$ is strictly decreasing then
$L_{w,1}$ has the Kadec-Klee property, i.e., if $f_n \rightarrow f$
weakly and $\|f_n\|_{w,1}\rightarrow \|f\|_{w,1}$ then
$\|f-f_n\|_{w,1}\rightarrow 0$. It is well-known, on the other hand,
that $L_1$ does not have this property. \endremark \head 5. Positive
Definite Functions \endhead

We would like to mention how the problems considered in this paper are
connected to positive definite functions.  The question we wish to
discuss is as follows: for which $q\in (0,2]$ and for which
 $a_1,\dots,a_n\geq 0$ is the function $\exp(-a_1 (x_1^{*})^q -\dots-
a_n (x_n^{*})^q)$ positive definite on $\Bbb R^n ?$ (Note that, for
every $q>2$, this function is not positive definite because its
one-dimensional restriction $\exp(-|t|^q)$ ($t\in \Bbb R$) is not
positive definite for these values of $q$.)

If $a_1=1$ and $a_2=\dots=a_n=0$ we arrive at the following question:
for which $q$ is the function $\exp(-\max^q(|x_1|,\dots,|x_n|))$
positive definite? This is exactly the problem posed in 1938 by I. J.
Schoenberg \cite{18} and solved by J. Misiewicz \cite{16} in 1989. The
answer is that, for every $n\geq 3$ and  $q>0$, the function is not
positive definite, and, for $n=2$, the function is positive definite if
and only if $q\leq 1.$ The problems of Schoenberg's type are important
in the study of isotropic and stable random vectors (for details, see
\cite{10}).

The connection between positive definite functions and isometric
embeddings into $L_q$ was discovered in 1966 by J. Bretagnolle, D.
Dacunha-Castelle and J. L. Krivine \cite{2}: for $q\in (0,2]$, a Banach
space $(E,\|\cdot\|)$ is isometric to a subspace of $L_q$ if and only
if the function $\exp(-\|x\|^q)$ is positive definite.

Combining this result with Theorem 1 we get an answer to the question
raised above.

\proclaim{Proposition 3} The function   $\exp(-a_1 (x_1^{*})^q -\dots -
a_n (x_n^{*})^q)$  is positive definite if and only if :  (i) $q\in
(1,2]$ and $a_1=\dots=a_n$, or \newline (ii) $q\leq 1$ and the numbers
$a_k$ form a non-increasing arithmetic progression. (This includes the
case $n=2$ where the progression consists of two numbers only.)
\endproclaim

Since we no longer assume in the formulation of Proposition 3 that
$a_1\geq\dots\geq a_n$ (which was necessary in the definition of the
Lorentz spaces) we have to make a remark concerning the proof.
Formally, without the condition $a_1\geq\dots\geq a_n$ the function
$u(x)=(a_1 (x_1^{*})^q +\dots+ a_n (x_n^{*})^q)^{1/q}$ may not be a
norm and so we cannot apply the result from \cite{2} directly. We can,
however, use a generalization of this result from \cite{17, p.290}
stating that, for any continuous one-homogeneous non-negative function
$u$ on $\Bbb R^n$ which vanishes only at the origin, if
$\exp(-{u(x)}^q)$ is positive definite, then $u$ is the norm (or
$q$-norm if $q<1$) of a subspace of $L_q.$ Besides, the proof of
Theorem~1 does not depend on the condition  $a_1\geq\dots\geq a_n$,
and, in fact, in the proof of Theorem 1 this condition follows from the
existence of an isometric embedding into $L_q.$

\enddocument

\bye

\bye